\documentclass{article}
\usepackage[latin1]{inputenc}
\usepackage{amsmath}
\usepackage{amsfonts}
\usepackage{amssymb}
\author{Vehbi Emrah Paksoy}
\usepackage{latexsym}
\usepackage{amsfonts}
\usepackage{amsmath}
\usepackage{xypic}
\newtheorem{df}[subsection]{\textbf{Definition }}
\newtheorem{prop}[subsection]{\textbf{Proposition }}
\newtheorem{lem}[subsection]{\textbf{Lemma }}
\newtheorem{thm}[subsection]{\textbf{Theorem }}
\newtheorem{cor}[subsection]{\textbf{Corollary }}

\title{On the integrability of spatial polygons in 3 dimensional Minkowsi Spaces}
\date{}
\begin{document}
\maketitle

\begin{abstract}
In this paper we show that the space of spatial polygons in 3 dimensional Minkowski space $\mathbb{M}^{3}$  is a K\"ahler manifold. We will give explicit description to the tangent space and constructing an almost complex structure we will show the integrability using Newlander-Nierenberg theorem. This approach certainly has lots of computational advantages.
\end{abstract}

Space of spatial polygones in Euclidean spaces are well understood after the works of \cite{HK,KM1,KLY1}. Similar constructions can be extended to the spatial polygones in semi-Riemannian spaces . We can define the Lorentz metric on $\mathbb{R}^{n}$ ,$n>1$ as $\displaystyle (u,v)=u_{n}v_{n}-\sum_{i=1}^{n-1}u_{i}v_{i}$ where $u=(u_{1},\ldots ,u_{n}), v=(v_{1},\ldots ,v_{n}) \in \mathbb{R}^{n}$. The vector space ($\mathbb{R}^{n},(,) $ ) is called Minkowski space and denoted by $\mathbb{M}^{n}$.
For those who interested in the algebraic and geometric properties of Minkowski spaces, \cite{Iv} and \cite{Ra} would be  good sources. Next Consider $\mathfrak{sl}_{2}(\mathbb{R})=\lbrace A \in M^{2\times 2}| trA=0\rbrace$. The form defined as $(A,B)=-\frac{1}{2}tr{AB}$ where $A,B \in \mathfrak{sl}_{2}(\mathbb{R})$ is a metric on $\mathfrak{sl}_{2}(\mathbb{R}).$ We know that the matrices
\begin{displaymath}
\mathbf{e_{1}}= \left (\begin{array}{ccc}
0 & 1 \\
1 & 0 
\end{array} \right ), 
\mathbf{e_{2}}= \left (\begin{array}{ccc}
1 & 0 \\
0 & -1
\end{array} \right ),
\mathbf{e_{3}}= \left ( \begin{array}{ccc}
0 & -1 \\
1 & 0 \\
\end{array} \right ).
\end{displaymath}
form a basis for $\mathfrak{sl}_{2}(\mathbb{R})$. Note that 
$$(e_{1},e_{2})=(e_{1},e_{3})=(e_{2},e_{3})=0$$ \\
moreover $(e_{1},e_{1})=(e_{2},e_{2})=-1$ and $(e_{3},e_{3})=1$. We conclude that $(\mathfrak{sl}_{2}(\mathbb{R}),(,))$ is of Sylvester type (-2,1). Let us define the map
\begin{eqnarray}
f: \mathbb{M}^{3} &\longrightarrow& \mathfrak{sl}_{2}(\mathbb{R}) \nonumber \\  
u=(x,y,z) &\mapsto& A \nonumber
\end{eqnarray}
Where $A$ is the matrix given by;
\begin{displaymath}
A=
\left (\begin{array}{ccc}
x & y+z \\
y-z & -x 
\end{array} \right ).
\end{displaymath}
$f$ is an isomorphism and $(fu,fv)=(u,v)$ for $u,v\in \mathbb{M}^{3}$. So $\mathbb{M}^{3} \simeq \mathfrak{sl}_{2}(\mathbb{R})$ with \linebreak $(A,B)=-\frac{1}{2}tr(AB)$ i.e, an isometry.

From now on, we will only consider $\mathbb{M}^{3}$ or equivalently $\mathfrak{sl}_{2}(\mathbb{R})$ with metric $(A,B)=-\frac{1}{2}tr(AB)$. On $\mathbb{M}^{3}$ we can define the product $[u,v]=J(u \times v)$ where
\begin{displaymath}
J=
\left ( \begin{array}{ccc}
1 & 0 & 0 \\
0 & -1 & 0 \\
0 & 0 & -1
\end{array} \right ).
\end{displaymath}
and `$\times$` denotes the usual vector product in $\mathbb{R}^{3}$. It will be useful to explore some properties of this product. The proof of the following proposition is obtained using definition.
\begin{prop} The product \textrm{[,]} on $\mathbb{M}^{3} $ posesses the following properties;
\begin{eqnarray}
i)&& [u,v] \mbox{ is orthagonal both } u \mbox{ and } v \mbox{ w.r.t Minkowski metric.} \nonumber \\
ii)&& [u,v]=-[v,u] \nonumber \\
iii)&& ([u,v],w)=(u,[v,w]) \nonumber \\
iv)&& [u,[v,w]]=(u,w)v-(u,v)w \nonumber \\
v)&& [u,v]=0 \Leftrightarrow u=\lambda v \ \lambda \in \mathbb{R} \nonumber
\end{eqnarray}
where $u,v,w \in \mathbb{M}^{3}$.
\end{prop}

Moreover $\mathbb{M}^{3}$ with the product [,] becomes a Lie algebra. Now recall that $\mathbb{M}^{3} \simeq \mathfrak{sl}_{2}(\mathbb{R})$ with given metric. We will define a product [,] on $\mathfrak{sl}_{2}(\mathbb{R})$ to be
$$[A,B]=f[f^{-1}(A),f^{-1}(B)]$$
where $f^{-1}(A)=(x,\frac{y+z}{2},\frac{y-z}{2})$ for
\begin{displaymath}
A=\left ( \begin{array}{cc}
x & y\\
z & -x
\end{array} \right). \in \mathfrak{sl}_{2}(\mathbb{R}).
\end{displaymath}
Note that 
$$([A,B],A)=(f[f^{-1}(A),f^{-1}(B)],A)=(f[f^{-1}(A),f^{-1}(B)],f(f^{-1}(A))).$$\\
Since $f$ is an isometry we have 
$$([A,B],A)=([f^{-1}(A),f^{-1}(B)],f^{-1}(A))=0$$\\
and similarly $([A,B],B)=0$.
\begin{prop} Product [,] on $\mathfrak{sl}_{2}(\mathbb{R})$ with given metric has the following properties;
\begin{eqnarray}
i&)& [A,B] \mbox{ is orthagonal to both } A \mbox{ and }B,\nonumber \\
ii&)& [A,B]=-[B,A],\nonumber \\
iii&)&[A,B]=0 \Leftrightarrow A=\lambda B,\nonumber \\
iv&)& [A,[B,C]]=(A,C)B-(A,B)C,\nonumber \\
v&)& (A,[B,C])=([A,B],C).\nonumber
\end{eqnarray}
where $A,B,C\in \mathfrak{sl}_{2}(\mathbb{R})$ and $\lambda \in \mathbb{R}.$
\end{prop}
{\bf Proof:} We prove the item $i)$ just before the proposition. All of the other items listed above can be checked directly using the properties of [,] in $\mathbb{M}^{3}.\quad \Box$\\
As a result of the previous proposition we can say that the bracket on $\mathbb{M}^{3}$ and the product on $\mathfrak{sl}_{2}(\mathbb{R})$ have the same properties.Observe that for $A\in \mathfrak{sl}_{2}(\mathbb{R})$ with $(A,A)=1$ and any $B\in \mathfrak{sl}_{2}(\mathbb{R})$ such that $(A,B)=0$ we have \linebreak $[A,[A,B]]=(A,B)C-(A,A)B=-B$. Thus the operator 
\begin{eqnarray}
I_{A}:<A>^{\bot} &\longrightarrow& <A>^{\bot} \nonumber \\
B&\mapsto&[A,B] \nonumber
\end{eqnarray}
defines an almost complex structure on $<A>^{\bot}$, orthagonal complement of A in $\mathfrak{sl}_{2}(\mathbb{R})$.

\section{Geometry of the Moduli}

Let $\widetilde{\mathcal{M}}$ be the finite set of $p_{\alpha} \in \mathfrak{sl}_{2}(\mathbb{R})$ such that 
$\sum_{\alpha} p_{\alpha}=0 \mbox{ and }\linebreak (p_{\alpha},p_{\alpha})=m_{\alpha}^{2} ,m_{\alpha}>0$. Hence $p_{\alpha}$ is time-like for all $\alpha$.
We set $\mathcal{M}=\widetilde{\mathcal{M}}$ modulo adjoint action of $SL_{2}(\mathbb{R})=\widetilde{\mathcal{M}}/SL_{2}(\mathbb{R})$. This is the space of polygons in $\mathbb{M}^{3}$ with sides $p_{\alpha}$ and side lengths $m_{\alpha}$.Then the tangent space $\mathcal{T}(P)$ at point $P=\lbrace p_{\alpha}\rbrace$ consist of vectors $q=\lbrace q_{\alpha} \rbrace$ such that  
\begin{eqnarray}
\mathit{i})&&  (q_{\alpha}, p_{\alpha}) =0 \nonumber \\
\mathit{ii})&&  \sum_{\alpha} q_{\alpha} =0 \nonumber \\
\mathit{iii})&&  \mbox{Two systems } \lbrace q_{\alpha}\rbrace \mbox{ and } \lbrace w_{\alpha} \rbrace \mbox { represents the same tangent vector if } \nonumber \\
&&\exists x \in \mathfrak{sl}_{2}(\mathbb{R}) \mbox{ such that } w_{\alpha}=q_{\alpha}+[x,p_{\alpha}].\nonumber
\end{eqnarray}
Third condition defining the equivalence of the vectors can be interpreted as the infinitelstimal motion of the polygon as a whole.
\begin{prop}In each class ef equivalence given in $iii$ , there exist a unique representative $\lbrace q_{\alpha} \rbrace$ such that;
$$\sum_{\alpha} \frac{[q_{\alpha},p_{\alpha}]}{m_{\alpha}}=0 \qquad (*)$$.
\end{prop}
{\bf Proof:} To show the existance, consider $q={q_{\alpha}}$, by i) we know $q_{\alpha}\bot p_{\alpha}$ i.e, $q_{\alpha}$ is space-like vector for all $\alpha.\Rightarrow (q_{\alpha},q_{\alpha})<0$. So there exist a representative $\{ q_{\alpha} \}$ such that 
$$-\sum_{\alpha} \frac{(q_{\alpha},q_{\alpha})}{m_{\alpha}}=min. $$\\
Hence by extrema condition we have $\sum_{\alpha}
\frac{[q_{\alpha},p_{\alpha}]}{m_{\alpha}}=0.$\\
For uniqueness, let $\widetilde{q_{\alpha}}=q_{\alpha}+[x,p_{\alpha}]$ for some $x\in \mathfrak{sl}_{2}(\mathbb{R})$. Note that  
$$\sum_{\alpha} 
\frac{[\widetilde{q_{\alpha}},p_{\alpha}]}{m_{\alpha}}=0.$$\\
Then
$$\sum_{\alpha} \frac{[[x,p_{\alpha}],p_{\alpha}]}{m_{\alpha}}=0 \Rightarrow \sum_{\alpha} \frac{([[x,p_{\alpha}],p_{\alpha}],x)}{m_{\alpha}}=0.\mbox{ Hence } \sum_{\alpha} \frac{([x,p_{\alpha}],[x,p_{\alpha}])}{m_{\alpha}}=0.$$\\
Since $p_{\alpha}$ and $[x,p_{\alpha}]$ are orthagonal, $[x,p_{\alpha}]$ is space-like. So \linebreak $([x,p_{\alpha}],[x,p_{\alpha}])<0 \Rightarrow [x,p_{\alpha}]=0$ for all $\alpha$. Hence $x=\lambda p_{\alpha}, \lambda \in \mathbb{R}$ for all $\alpha$. For non-collinear $p_{\alpha}$`s , this is the case if $x=0$. Therefore $\widetilde{q_{\alpha}}=q_{\alpha}. \qquad \Box$
In the case of collinear vectors, the gauge representative $\{q_{\alpha}\}$ is unique but the tangent space  has dimension greater that that of $\mathcal{M}$. This means that the point is singular.\\
The uniqueness allow us to define the operator 
\begin{eqnarray}
I:\mathcal{T}(P) &\longrightarrow& \mathcal{T}(P) \nonumber \\
\{q_{\alpha} \} &\mapsto& \{ \frac{[q_{\alpha},p_{\alpha}]}{m_{\alpha}} \} \nonumber
\end{eqnarray}
where $q_{\alpha}$ satisfies the calibration condition $\sum_{\alpha} 
\frac{[q_{\alpha},p_{\alpha}]}{m_{\alpha}}=0$. Note that \linebreak $I^{2}q_{\alpha}=I(Iq_{\alpha})=-q_{\alpha}$. Therefore we obtain a complex structure on tangent space $\mathcal{T}(P)$.
\begin{lem} The 2-form $\omega$ on $\mathcal{T}(P)$ given by
$$\omega(q,q^{\prime})=\sum_{\alpha} 
\frac{([q_{\alpha},q^{\prime}_{\alpha}],p_{\alpha})}{m_{\alpha}^{2}}$$ \\
where $q=\{ q_{\alpha} \},q^{\prime}=\{ q^{\prime}_{\alpha} \} \in \mathcal{T}(P)$ is a symplectic form and in addition invariant under gauge transformation $q_{\alpha} \mapsto q_{\alpha}+[a,p_{\alpha}],a\in \mathfrak{sl}_{2}(\mathbb{R})=\mathbb{M}^{3}$.
\end{lem}
{\bf Proof:} Let us first prove the invariance. Consider
\begin{eqnarray}
&&\hspace{-1cm} \sum_{\alpha}
\frac{([q_{\alpha}+[a,p_{\alpha}],q^{\prime}_{\alpha}+[b,p_{\alpha}]],p_{\alpha})}{m_{\alpha}^{2}}=\sum_{\alpha}
\frac{(q_{\alpha}+[a,p_{\alpha}],[q^{\prime}_{\alpha}+[b,p_{\alpha}],p_{\alpha}])}{m_{\alpha}^{2}}=\nonumber \\
=&&-\sum_{\alpha}
\frac{(q_{\alpha}+[a,p_{\alpha}],m_{\alpha}^{2}b-(p_{\alpha},b)p_{\alpha}+[p_{\alpha},q^{\prime}_{\alpha}])}{m_{\alpha}^{2}}= \nonumber\\
=&&\sum_{\alpha}
\frac{([q_{\alpha},q^{\prime}_{\alpha}],p_{\alpha})}{m_{\alpha}^{2}}=\omega(q,q^{\prime}). \nonumber
\end{eqnarray}
Hence, $\omega$ is invariant under gauge transformation. Note that 
$$\omega(q,q^{\prime})=\sum_{\alpha} 
\frac{([q_{\alpha},q^{\prime}_{\alpha}],p_{\alpha})}{m_{\alpha}^{2}}=-\omega(q^{\prime},q)$$
So, $\omega$ is alternating.\\
Now assume $\omega(q,q^{\prime})=0, \forall q^{\prime} \in \mathcal{T}(P)$. Then set $q^{\prime}=Iq \in \mathcal{T}(P)$ we get
\begin{eqnarray}
0=&&\sum_{\alpha} 
\frac{([q_{\alpha},Iq_{\alpha}],p_{\alpha})}{m_{\alpha}^{2}}=\sum_{\alpha}
\frac{([q_{\alpha},[q_{\alpha},p_{\alpha}]],p_{\alpha})}{m_{\alpha}^{3}}=\nonumber \\
=&&-\sum_{\alpha}\frac{(q_{\alpha},q_{\alpha})(p_{\alpha},p_{\alpha})}{m_{\alpha}^{3}}=-\sum_{\alpha} \frac{(q_{\alpha},q_{\alpha})}{m_{\alpha}} \nonumber
\end{eqnarray}
Note that $q_{\alpha}$`s are space-like vectors. So, $(q_{\alpha},q_{\alpha})<0$. In this case above equality holds iff $q_{\alpha}=0, \forall \alpha$.This means that $\omega$ is non-degenerate. \\
It remains to show that $\omega$ is closed. Observe that it is closed on `mass surfaces' $(p_{\alpha},p_{\alpha})=m_{\alpha}^{2}$ and by the invariance, it is also closed on the factor.\qquad $\Box$\\
After defining symplectic form $\omega$  on $\mathcal{M}$ it is convenient to construct a Riemannian metric using $\omega$. Let`s define;
$$g(q,q^{\prime})=-\omega(Iq,q^{\prime})$$\\
Then $g$ is non-degenerate symmetric form and for the vectors satisfying calibration condition in proposition(3.17), it can be written as
$$g(q,q^{\prime})=-\sum_{\alpha} 
\frac{(q_{\alpha},q^{\prime}_{\alpha})}{m_{\alpha}}$$ \\
and $g(q,q)>0$ since $(q_{\alpha},q_{\alpha})<0$.

\section{Integrability}

In previous section we defined the almost complex structure $I$ on $\mathcal{M}$ and tangent space at a point $P \in \mathcal{M}$. We may consider the tangent space $\mathcal{T}(P)$ as a subset in $(\mathbb{M}^{3})^{n}$ where $\alpha \in \{ 1,\ldots ,n\}.$ Consider the 2-form defined by
$$(x,y)=\sum_{\alpha} \frac{(x_{\alpha},y_{\alpha})}{m_{\alpha}}$$ \\
where $x=\{ x_{\alpha} \},y=\{ y_{\alpha} \} \in (\mathbb{M}^{3})^{n}$ and $(x_{\alpha},y_{\alpha})$ is the Minkowski metric on $\mathbb{M}^{3}$.
Let $\mathcal{N}_{\mbox{P}}$ be the orthagonal complement to $\mathcal{T}(P)$ via the form defined above.\linebreak
So 
\begin{displaymath}
(\mathbb{M}^{3})^{n}=\mathcal{T} \mbox{(P)} \oplus \mathcal{N}_{\mbox{P}}
\end{displaymath}
\begin{prop} The operator $L:\mathbb{M}^{3} \longrightarrow 
\mathbb{M}^{3}$ is given by;
$$L(\xi)=\sum_{\alpha} \frac{[p_{\alpha},[\xi,p_{\alpha}]]}{m_{\alpha}}$$\\ 
is self-adjoint with respect to Minkowski metric.
\end{prop}
{\bf Proof:} Consider
\begin{eqnarray}
(L(\xi),\eta)=\sum_{\alpha} (\frac{[p_{\alpha},[\xi,p_{\alpha}]]}{m_{\alpha}},\eta)
&=& \sum_{\alpha} (\frac{[p_{\alpha},[\eta,p_{\alpha}]]}{m_{\alpha}},\xi) \nonumber\\
&=& (\sum_{\alpha} \frac{[p_{\alpha},[\eta,p_{\alpha}]]}{m_{\alpha}},\xi) \nonumber \\
&=& (L(\eta),\xi) \qquad \Box \nonumber
\end{eqnarray}
The following lemma gives a explicit projection to the normal $\mathcal{N}_{P}$.
\begin{lem} Let $\pi:(\mathbb{M}^{3})^{n} \longrightarrow 
\mathcal{N}_{P}$ given by 
$$(\pi x)_{\alpha}=\frac{(x_{\alpha},p_{\alpha})}{m_{\alpha}^{2}} 
p_{\alpha}+\frac{[p_{\alpha},[\xi_{x},p_{\alpha}]]}{m_{\alpha}}+
[w_{x},p_{\alpha}]$$
where $\xi_{x}, w_{x} \in \mathbb{M}^{3}$ is defined uniquely such that
\begin{eqnarray}
L(\xi_{x})&=& \sum_{\alpha} \frac{[p_{\alpha},[\xi_{x},p_{\alpha}]]}{m_{\alpha}}= 
\sum_{\alpha} \frac{[p_{\alpha},[x_{\alpha},p_{\alpha}]]}{m_{\alpha}^{2}} \nonumber \\ 
L(w_{x})&=& \sum_{\alpha} \frac{[p_{\alpha},[w_{x},p_{\alpha}]]}{m_{\alpha}}= 
\sum_{\alpha} \frac{[p_{\alpha},x_{\alpha}]}{m_{\alpha}} \nonumber
\end{eqnarray}
$\alpha \in \lbrace 1,\ldots,\mbox{n} \rbrace  \mbox{ and }\pi \mbox{ is 
self-adjoint.}$
\end{lem}
{\bf Proof: } We must show that $(\pi x,u)= 0, \forall u \in \mathcal{T}(P).$ 
Let \linebreak $u=(u_{1},\ldots,u_{n})\in \mathcal{T}(P)$
\begin{eqnarray}
(\pi x,u)&=& \sum_{\alpha}\Big ( (\frac{(x_{\alpha},p_{\alpha})}{m_{\alpha}^{3}} 
p_{\alpha},u_{\alpha})+
(\frac{[p_{\alpha},[\xi_{x},p_{\alpha}]]}{m_{\alpha}^{2}},u_{\alpha})+
(\frac{[w_{x},p_{\alpha}]}{m_{\alpha}},u_{\alpha}) \Big ) \nonumber \\
&=& -\sum_{\alpha} 
\frac{(p_{\alpha},[u_{\alpha},[\xi_{x},p_{\alpha}]])}{m_{\alpha}^{2}}+ 
\sum_{\alpha} \frac{(w_{x},[p_{\alpha},u_{\alpha}])}{m_{\alpha}} = 0.\nonumber
\end{eqnarray}
Since $(u_{\alpha},p_{\alpha})=\sum_{\alpha} 
\frac{[p_{\alpha},u_{\alpha}]}{m_{\alpha}}=\sum_{\alpha} u_{\alpha}=0.$
Consider $x,y \in (\mathbb{M}^{3})^{n}.$ We know that 
$$(y,\pi x)=\sum_{\alpha} \Big ( 
(y_{\alpha} ,\frac{(p_{\alpha},x_{\alpha})}{m_{\alpha}^{3}}p_{\alpha})+(y_{\alpha},
\frac{[p_{\alpha},[\xi_{x},p_{\alpha}]]}{m_{\alpha}^{2}})+
(y_{\alpha},\frac{[w_{x},p_{\alpha}]}{m_{\alpha}}\mbox{) } \Big )$$\\
So
$$(y,\pi x)=\sum_{\alpha} 
\Big ( \frac{(y_{\alpha},p_{\alpha})(p_{\alpha},x_{\alpha})}{m_{\alpha}^{3}}+
(L(\xi_{y}),\xi_{x})+(L(w_{y}),w_{x}) \Big ) \hspace{1.7cm} \mbox{(I)}$$
 and
\begin{eqnarray}
(\pi y,x)&=& \sum_{\alpha} 
(x_{\alpha},\frac{(p_{\alpha},y_{\alpha})}{m_{\alpha}^{3}}p_{\alpha})  
+\sum_{\alpha}(x_{\alpha},
\frac{[p_{\alpha},[\xi_{y},p_{\alpha}]]}{m_{\alpha}^{2}})+
\sum_{\alpha} (x_{\alpha},\frac{[w_{y},p_{\alpha}]}{m_{\alpha}}) \nonumber \\
&=& \sum_{\alpha} 
\Big ( \frac{(x_{\alpha},p_{\alpha})(y_{\alpha},p_{\alpha})}{m_{\alpha}^{3}}
+(L(\xi_{x}),\xi_{y})+(L(w_{x}),w_{y})\Big ) \hspace{1.5cm} \mbox{(II)} \nonumber
\end{eqnarray}
Comparing I and II we obtain $(\pi x,y)=(x,\pi y) \qquad \Box$ \\
We know that edges of a polygon from moduli space are time-like vectors 
and since components of tangents are Lorentz orthagonal to them, these 
components are space-like vectors. Hence we have a Riemannian metric on 
tangent space $\mathcal{T} \mbox{(P) given by; }$
$$g(u,v)=- \sum_{\alpha} \frac{(u_{\alpha},v_{\alpha})}{m_{\alpha}}$$
Symmetricity and bilinearity of this form easily obtained from 
properties of Lorentzian metric. We know that components of each tangent 
vector $u \in \mathcal{T} \mbox{(P)}$ is space-like. So $(u_{\alpha},u_{\alpha})
<0 \Rightarrow g(u,u)>0$. Thus, $g(u,v)$ is positive and so is a Riemannian metric on tangent space at point P.

We can define covariant derivative of a vector field 
$q$ along a curve $P=P(t)$ in $\mathcal{M}$ by;
$$\nabla_{t}q=\mbox{ projection }\Big ( \frac{dP}{dt} \in 
(\mathbb{M}^{3})^{n} \Big ) \mbox{ in tangent space of } \mathcal{M}$$\\
Therefore if $x,y$  are vector fields then at a point $P$
Lie bracket $ [x,y]=\nabla_{x}y-\nabla_{y}x \in \mathcal{T} 
\mbox{(P)}$
\begin{prop} Let $x, y$ be vector fields. At a point $P\in \mathcal{M}$ we have the following
$$\{(\nabla_{x}y)^{\alpha}\}= 
\{\partial_{x}y^{\alpha}-\pi(\partial_{x}y)_{\alpha} \} \in \mathcal{T}(P)$$. Here $\partial$ is the covariant derivative on $(\mathbb{M}^{3})^{n}$ and $\partial_{x}y^{\alpha}$ denotes $\alpha$ component of $\partial_{x}y$ for $\alpha=1,\ldots ,n$.
\end{prop}
{\bf Proof:} Note that for $x, y$ vector fields we have ;
$$(\nabla_{x}y)^{\alpha} 
=\partial_{x}y^{\alpha}-
\frac{(p_{\alpha},\partial_{x}y^{\alpha})p_{\alpha}}{m_{\alpha}^{2}}-
\frac{[p_{\alpha},[\xi_{\partial_{x}y},p_{\alpha}]]}{m_{\alpha}}-
[w_{\partial_{x}y},p_{\alpha}].$$ \\
We shall show that $\lbrace (\nabla_{x}y)^{\alpha}\rbrace  \mbox{ satisfies 
\textit {i-iii}}.$ 
\begin{eqnarray}
i)&&\sum_{\alpha} (\nabla_{x}y)^{\alpha}= \sum_{\alpha} \partial_{x}y^{\alpha} -\sum_{\alpha}
\frac{(p_{\alpha},\partial_{x}y^{\alpha})p_{\alpha}}{m_{\alpha}^{2}}-
\sum_{\alpha} 
\frac{[p_{\alpha},[\xi_{\partial_{x}y},p_{\alpha}]]}{m_{\alpha}}-\sum_{\alpha} [w_{\partial_{x}y},p_{\alpha}]= \nonumber \\
&& -\sum_{\alpha} 
\frac{(p_{\alpha},\partial_{x}y^{\alpha})p_{\alpha}}{m_{\alpha}^{2}}-
\sum_{\alpha} 
\frac{[p_{\alpha},[\partial_{x}y^{\alpha},p_{\alpha}]]}{m_{\alpha}^{2}}=0.\nonumber \\
ii)&&\mbox{We must show that }((\nabla_{x}y)^{\alpha},p_{\alpha})=0.\mbox{ Note that}\nonumber \\
&&(\partial_{x}y^{\alpha}-
\frac{(\partial_{x}y^{\alpha},p_{\alpha})p_{\alpha}}{m_{\alpha}^{2}}-
\frac{[p_{\alpha},[\xi_{\partial_{x}y},p_{\alpha}]]}{m_{\alpha}}-
[w_{\partial_{x}y},p_{\alpha}],p_{\alpha})=\nonumber \\ 
&&(\partial_{x}y^{\alpha},p_{\alpha})-(\partial_{x}y^{\alpha},p_{\alpha})=0.\nonumber
\end{eqnarray}
\textit{iii)} \hspace{0.5cm} Consider 
\begin{eqnarray}
\sum_{\alpha} \frac{[(\nabla_{x}y)^{\alpha},p_{\alpha}]}{m_{\alpha}}  &=& \sum_{\alpha} \Big ( \frac{[\partial_{x}y^{\alpha},p_{\alpha}]}{m_{\alpha}}+
[\xi_{\partial_{x}y},p_{\alpha}]-
\frac{[[w_{\partial_{x}y},p_{\alpha}],p_{\alpha}]}{m_{\alpha}} \Big )  \nonumber \\
&=& \sum_{\alpha} \frac{[\partial_{x}y^{\alpha},p_{\alpha}]}{m_{\alpha}}+\sum_{\alpha} 
\frac{[p_{\alpha},[w_{\partial_{x}y},p_{\alpha}]]}{m_{\alpha}}=0 \nonumber 
\hspace{1cm}  \Box 
\end{eqnarray}
Let us define
\begin{eqnarray}
\mu: \mathcal{T}\mathcal{M}\times \mathcal{T} \mathcal{M} &\rightarrow&  \mathbb{M}^{3}  \nonumber \\ 
\hspace{-3cm}
(x,y) &\mapsto& \mu(x,y)  \nonumber 
\end{eqnarray}
Where $\mu(x,y)$ is the vector satisfying 
$$\sum_{\alpha} 
\frac{[p_{\alpha},[\mu(x,y),p_{\alpha}]]}{m_{\alpha}}=\sum_{\alpha} 
\frac{[x_{\alpha},[y_{\alpha},p_{\alpha}]]}{m_{\alpha}^{2}}.$$
Note $\partial_{x}p_{\alpha}$ is the $\alpha$ component of $x(P)\in \mathcal{T}(P)$. So $\partial_{x}p_{\alpha}=x_{\alpha}$.
\begin{prop} For $x,y$ vector fields on $\mathcal{M}$, we have 
$$\xi_{\partial_{x}y}=-\mu(x,y) \mbox{ and } w_{\partial_{x}y}=\mu(Iy,x)$$
\end{prop}
{\bf Proof:}  
\begin{eqnarray}
\sum_{\alpha} 
\frac{[p_{\alpha},[\partial_{x}y^{\alpha},p_{\alpha}]]}{m_{\alpha}^{2}} &=&
 - \sum_{\alpha} 
\frac{(p_{\alpha},\partial_{x}y^{\alpha})p_{\alpha}}{m_{\alpha}^{2}}=
\sum_{\alpha} \frac{(x_{\alpha},y_{\alpha})p_{\alpha}}{m_{\alpha}^{2}} \nonumber \\
&=& -\sum_{\alpha} 
\frac{[x_{\alpha},[y_{\alpha},p_{\alpha}]]}{m_{\alpha}^{2}}=-\sum_{\alpha}
\frac{[p_{\alpha},[\mu (x,y),p_{\alpha}]]}{m_{\alpha}} \nonumber \\
&=& \sum_{\alpha} \frac{[p_{\alpha},[-\mu (x,y),p_{\alpha}]]}{m_{\alpha}} \nonumber
\end{eqnarray}

\noindent Since $\xi_{\partial_{x}y}$ is the unique vector 
satisfying  
$$\sum_{\alpha} 
\frac{[p_{\alpha},[\xi_{\partial_{x}y},p_{\alpha}]]}{m_{\alpha}}=\sum_{\alpha}
\frac{[p_{\alpha},[\partial_{x}y^{\alpha},p_{\alpha}]]}{m_{\alpha}^{2}},$$ 
 we have  $ \xi_{\partial_{x}y}=-\mu(x,y).$\\
Now consider
$$\sum_{\alpha} 
\frac{[p_{\alpha},[w_{\partial_{x}y},p_{\alpha}]]}{m_{\alpha}}=\sum_{\alpha}
\frac{[p_{\alpha},\partial_{x}y^{\alpha}]}{m_{\alpha}}$$
We know that $\sum_{\alpha} 
\frac{[p_{\alpha},y_{\alpha}]}{m_{\alpha}}=0.$ 
 Differentiatingthis in the direction of $x$ we get;
$$- \sum_{\alpha} \frac{[x_{\alpha},y_{\alpha}]}{m_{\alpha}}=\sum_{\alpha}
\frac{[p_{\alpha},\partial_{x}y^{\alpha}]}{m_{\alpha}}$$
so
$$\sum_{\alpha} 
\frac{[p_{\alpha},[w_{\partial_{x}y},p_{\alpha}]]}{m_{\alpha}}=
-\sum_{\alpha} \frac{[x_{\alpha},y_{\alpha}]}{m_{\alpha}}$$\\ 
We have $[p_{\alpha},[x_{\alpha},y_{\alpha}]]=0 
\Rightarrow [x_{\alpha},y_{\alpha}]=\lambda_{\alpha}p_{\alpha}.$ Hence
$(p_{\alpha},[x_{\alpha},y_{\alpha}])=\lambda_{\alpha}m_{\alpha}^{2} \mbox{ and } 
\lambda_{\alpha}=\displaystyle
\frac{(p_{\alpha},[x_{\alpha},y_{\alpha}])}{m_{\alpha}^{2}}
$\\
So we get
$$-\sum_{\alpha} 
\frac{[x_{\alpha},y_{\alpha}]}{m_{\alpha}}=-\sum_{\alpha} 
\frac{(p_{\alpha},[x_{\alpha},y_{\alpha}])p_{\alpha}}{m_{\alpha}^{3}} 
$$\\
Now observe that
\begin{eqnarray}
\sum_{\alpha} 
\frac{[p_{\alpha},[\mu(Iy,x),p_{\alpha}]]}{m_{\alpha}}=&& \sum_{\alpha}
\frac{[Iy_{\alpha},[x_{\alpha},p_{\alpha}]]}{m_{\alpha}^{2}}=
 \sum_{\alpha} 
\frac{[[y_{\alpha},p_{\alpha}],[x_{\alpha},p_{\alpha}]]}{m_{\alpha}^{3}}= 
\nonumber \\
=&&-\sum_{\alpha} 
\frac{([y_{\alpha},p_{\alpha}],x_{\alpha})p_{\alpha}}{m_{\alpha}^{3}}= 
 -\sum_{\alpha} 
\frac{(p_{\alpha},[x_{\alpha},y_{\alpha}])p_{\alpha}}{m_{\alpha}^{3}}=\nonumber \\
=&&-\sum_{\alpha} \frac{[x_{\alpha},y_{\alpha}]}{m_{\alpha}} \nonumber
\end{eqnarray}
Therefore 
$$\sum_{\alpha} 
\frac{[p_{\alpha},[\mu(Iy,x),p_{\alpha}]]}{m_{\alpha}}=-\sum_{\alpha} 
\frac{[x_{\alpha},y_{\alpha}]}{m_{\alpha}}=\sum_{\alpha} 
\frac{[p_{\alpha},\partial_{x}y^{\alpha}]}{m_{\alpha}}$$
but $w_{\partial_{x}y}$ is the unique vector satisfying 
$$\sum_{\alpha} 
\frac{[p_{\alpha},[w_{\partial_{x}y},p_{\alpha}]]}{m_{\alpha}}=
\sum_{\alpha} \frac{[p_{\alpha},\partial_{x}y^{\alpha}]}{m_{\alpha}}
\Rightarrow w_{\partial_{x}y}=\mu(Iy,x) \hspace{0.5cm} \Box $$

\noindent Using this proposition we may write;\\
$$(\nabla_{x}y)^{\alpha}=\partial_{x}y^{\alpha}-
\frac{(\partial_{x}y^{\alpha},p_{\alpha})p_{\alpha}}{m_{\alpha}^{2}}+
\frac{[p_{\alpha},[\mu(x,y),p_{\alpha}]]}{m_{\alpha}}-[\mu(Iy,x),p_{\alpha}]$$
The map $\mu(x,y) \mbox{ restricted on } \mathcal{T}(\mbox{P}) \mbox{ 
where P} \in \mathcal{M}$  has a lot of useful properties. The 
following lemma is devoted to write these properties explicitely.
\begin{lem}
$\mbox{For the map }\mu(x,y) \mbox{ restricted on }\mathcal{T}(\mbox{P}) 
\mbox{ we have the following;}$
\begin{eqnarray}
\mathit{i})&& \mu(x,y) \mbox{ is the unique vector satisfying } \nonumber \\
&& \hspace{2cm} \sum_{\alpha} 
\frac{[p_{\alpha},[\mu(x,y),p_{\alpha}]]}{m_{\alpha}}=\sum_{\alpha} 
\frac{[x_{\alpha},[y_{\alpha},p_{\alpha}]]}{m_{\alpha}^{2}}
\nonumber \\
\mathit{ii})&& \mu(x,y)=\mu(y,x) \nonumber \\
\mathit{iii})&& \mu \mbox{ is skew-adjoint, i.e, } \mu(Ix,y)=-\mu(x,Iy) \nonumber \\
\mathit{iv})&& \mu(Ix,Iy)= \mu(x,y) \nonumber
\end{eqnarray}
Where $x,y$ are vector fields in $\mathcal{M}.$
\end{lem}
{\bf Proof:}\textit{i}) Since $\mu(x,y)=-\xi_{\partial_{x}y}$ and $\xi_{\partial_{x}y}$ is unique we are done.\\
\textit{ii}) Consider
\begin{eqnarray}
\sum_{\alpha} \frac{[p_{\alpha},[\mu(y,x),p_{\alpha}]]}{m_{\alpha}}=&& \sum_{\alpha} \frac{[y_{\alpha},[x_{\alpha},y_{\alpha}]]}{m_{\alpha}^{2}}=-\sum_{\alpha}
\frac{(y_{\alpha},x_{\alpha})p_{\alpha}}{m_{\alpha}^{2}} \nonumber \\
\sum_{\alpha} \frac{[p_{\alpha},[\mu(x,y),p_{\alpha}]]}{m_{\alpha}}=&& \sum_{\alpha} \frac{[x_{\alpha},[y_{\alpha},p_{\alpha}]]}{m_{\alpha}^{2}}=-\sum_{\alpha}
\frac{(x_{\alpha},p_{\alpha})p_{\alpha}}{m_{\alpha}^{2}}. \nonumber
\end{eqnarray}
By uniqueness of $\mu(x,y)$ we have $\mu(x,y)=\mu(y,x)$\\
\textit{iii}) Note that
\begin{eqnarray}
\sum_{\alpha} \frac{[p_{\alpha},[\mu(Ix,y),p_{\alpha}]]}{m_{\alpha}}=&&\sum_{\alpha} \frac{[Ix_{\alpha},[y_{\alpha},p_{\alpha}]]}{m_{\alpha}^{2}}=-\sum_{\alpha}
\frac{([x_{\alpha},p_{\alpha}],y_{\alpha})p_{\alpha}}{m_{\alpha}^{3}} \nonumber \\
\sum_{\alpha} \frac{[p_{\alpha},[-\mu(x,Iy),p_{\alpha}]]}{m_{\alpha}}=&&-\sum_{\alpha} \frac{[x_{\alpha},[y_{\alpha},p_{\alpha}]]}{m_{\alpha}^{3}}=-\sum_{\alpha} \frac{([x_{\alpha},p_{\alpha}],y_{\alpha})p_{\alpha}}{m_{\alpha}^{3}}. \nonumber
\end{eqnarray}
Hence by uniqueness of $\mu(Ix,y)$ we have $\mu(Ix,y)=-\mu(x,Iy)$\\
\textit{iv}) Easily follows from \textit{iii}). \qquad $\Box$ \\
It is convienient to make an observation right now. Note that
$$\frac{[Iy_{\alpha},x_{\alpha}]}{m_{\alpha}}=\frac{(x_{\alpha},y_{\alpha})
p_{\alpha}}{m_{\alpha}^{2}}=
-\frac{(\partial_{x}y^{\alpha},p_{\alpha})p_{\alpha}}{m_{\alpha}^{2}}.$$
So we will get;
$$(\nabla_{x}y)^{\alpha}=\partial_{x}y^{\alpha}+
\frac{[Iy_{\alpha},x_{\alpha}]}{m_{\alpha}}-I[\mu(x,y),p_{\alpha}]-
[\mu(Iy,x),p_{\alpha}].$$ 

We have arrived the main result of this chapter. The following theorem is important in the sense that it reveals the geometrical structure of the variety of spatial polygons in $\mathbb{M}^{3}.$
\begin{thm}Almost complex structure I on $\mathcal{T}(\mbox{P})$ is integrable.
\end{thm}
{\bf Proof:} By theorem of Newlander-Nierenberg,\cite{NN}, it is enough to check Nijenhaus tensor $N_{I}(x,y)=0$ for all vector fields $x,y \mbox{ in } \mathcal{M}.$ Nijenhaus tensor is defined by ;
$$N_{I}=\frac{1}{2}([Ix,Iy]-[x,y]-I[Ix,y]-I[x,Iy]).$$\\
As a convention we write $[x,y]_{\alpha}=(\nabla_{x}y)^{\alpha}-(\nabla_{y}x)^{\alpha}$
\begin{eqnarray}
(\nabla_{Ix}Iy)^{\alpha}=&&\partial_{Ix}Iy^{\alpha}+
\frac{[IIy_{\alpha},Ix_{\alpha}]}{m_{\alpha}}-I[\mu(Iy,Ix),p_{\alpha}]-[\mu(IIy,Ix),p_{\alpha}] \nonumber \\
=&& \partial_{Ix}Iy^{\alpha}+\frac{[Ix_{\alpha},y_{\alpha}]}{m_{\alpha}}-
I[\mu(y,x),p_{\alpha}]+[\mu(y,Ix),p_{\alpha}] \nonumber \\
(\nabla_{Iy}Ix)^{\alpha}=&&\partial_{Iy}Ix^{\alpha}+
\frac{[IIx_{\alpha},Iy_{\alpha}]}{m_{\alpha}}-I[\mu(Ix,Iy),p_{\alpha}]-[\mu(IIx,Iy),p_{\alpha}] \nonumber \\
=&&\partial_{Iy}Ix^{\alpha}+\frac{[Iy_{\alpha},x_{\alpha}]}{m_{\alpha}}-I[\mu(x,y),p_{\alpha}]+[\mu(x,Iy),p_{\alpha}] \nonumber
\end{eqnarray} 
Note that
$$\frac{[Iy_{\alpha},x_{\alpha}]}{m_{\alpha}}=
\frac{[[y_{\alpha},p_{\alpha}],x_{\alpha}]}{m_{\alpha}^{2}}=
+\frac{(x_{\alpha},y_{\alpha})p_{\alpha}}{m_{\alpha}^{2}}$$ \\
and
$$\frac{[Ix_{\alpha},y_{\alpha}]}{m_{\alpha}}=
\frac{[[x_{\alpha},p_{\alpha}],y_{\alpha}]}{m_{\alpha}^{2}}=
+\frac{(y_{\alpha},x_{\alpha})p_{\alpha}}{m_{\alpha}^{2}}$$ \\
Comparing these two equations we see that 
$$\frac{[Iy_{\alpha},x_{\alpha}]}{m_{\alpha}}=
\frac{[Ix_{\alpha},y_{\alpha}]}{m_{\alpha}}$$ \\
So we get;
$$[Ix,Iy]_{\alpha}=\partial_{Ix}Iy^{\alpha}-\partial_{Iy}Ix^{\alpha}+2[\mu(y,Ix),p_{\alpha}]$$
\begin{eqnarray}
(\nabla_{x}y)^{\alpha}=&&\partial_{x}y^{\alpha}+
\frac{[Iy_{\alpha},x_{\alpha}]}{m_{\alpha}}-I[\mu(y,x),p_{\alpha}]-[\mu(Iy,x),p_{\alpha}] \nonumber \\
(\nabla_{y}x)^{\alpha}=&&\partial_{y}x^{\alpha}+
\frac{[Ix_{\alpha},y_{\alpha}]}{m_{\alpha}}-I[\mu(x,y),p_{\alpha}]-[\mu(Ix,y),p_{\alpha}] \nonumber
\end{eqnarray}
Therefore;
$$[x,y]_{\alpha}=\partial_{x}y^{\alpha}-\partial_{y}x^{\alpha}-2[\mu(Iy,x),p_{\alpha}]$$
\begin{eqnarray}
(\nabla_{Ix}y)^{\alpha}=&&\partial_{Ix}y^{\alpha}+
\frac{[Iy_{\alpha},Ix_{\alpha}]}{m_{\alpha}}-I[\mu(y,Ix),p_{\alpha}]-[\mu(Iy,Ix),p_{\alpha}]= \nonumber \\
=&&\partial_{Ix}y^{\alpha}+
\frac{[Iy_{\alpha},Ix_{\alpha}]}{m_{\alpha}}-I[\mu(y,Ix),p_{\alpha}]-[\mu(y,x),p_{\alpha}] \nonumber \\
(\nabla_{y}Ix)^{\alpha}=&& \partial_{y}Ix^{\alpha}+
\frac{[IIx_{\alpha},y_{\alpha}]}{m_{\alpha}}-I[\mu(Ix,y),p_{\alpha}]-[\mu(IIx,y),p_{\alpha}]= \nonumber \\
=&& \partial_{y}Ix^{\alpha}-
\frac{[x_{\alpha},y_{\alpha}]}{m_{\alpha}}-I[\mu(Ix,y),p_{\alpha}]+[\mu(x,y),p_{\alpha}] \nonumber 
\end{eqnarray}
Since $[Ix_{\alpha},y_{\alpha}]=[Iy_{\alpha},x_{\alpha}]$ we have;
$$[Iy_{\alpha},Ix_{\alpha}]=[IIx_{\alpha},y_{\alpha}]=-[x_{\alpha},y_{\alpha}]$$\\
so 
$$[Ix,y]_{\alpha}=\partial_{Ix}y^{\alpha}-\partial_{y}Ix^{\alpha}-2[\mu(x,y),p_{\alpha}]$$
\begin{eqnarray}
(\nabla_{x}Iy)^{\alpha}=&& \partial_{x}Iy^{\alpha}+
\frac{[IIy_{\alpha},x_{\alpha}]}{m_{\alpha}}-I[\mu(Iy,x),p_{\alpha}]-[\mu(IIy,x),p_{\alpha}] \nonumber \\
 =&& \partial_{x}Iy^{\alpha}-
\frac{[y_{\alpha},x_{\alpha}]}{m_{\alpha}}-I[\mu(Iy,x),p_{\alpha}]+[\mu(y,x),p_{\alpha}] \nonumber \\
(\nabla_{Iy}x)^{\alpha}=&&\partial_{Iy}x^{\alpha}+
\frac{[Ix_{\alpha},Iy_{\alpha}]}{m_{\alpha}}-I[\mu(x,Iy),p_{\alpha}]-[\mu(Ix,Iy),p_{\alpha}] \nonumber \\
=&&\partial_{Iy}x^{\alpha}-
\frac{[y_{\alpha},x_{\alpha}]}{m_{\alpha}}-I[\mu(x,Iy),p_{\alpha}]-[\mu(x,y),p_{\alpha}] \nonumber 
\end{eqnarray}
so
$$[x,Iy]_{\alpha}=\partial_{x}Iy^{\alpha}-\partial_{Iy}x^{\alpha}+2[\mu(x,y),p_{\alpha}]$$\\
With these tools we can calculate the Nijenhaus tensor $N_{I}(x,y)$. If we plug what we have found above we will obtain;
$$N_{I}(x,y)=\partial_{Ix}Iy^{\alpha}-\partial_{Iy}Ix^{\alpha}-\partial_{x}y^{\alpha}+\partial_{y}x^{\alpha}-I\partial_{Ix}y^{\alpha}+I\partial_{y}Ix^{\alpha}-
I\partial_{x}Iy^{\alpha}+I\partial_{Iy}x^{\alpha}$$\\
\textbf{Claim:} $N_{I}(x,y)=0$.\\
In order to show the claim we need to do some calculations. Observe that; 
$$[Iy_{\alpha},p_{\alpha}]= 
\frac{[[y_{\alpha},p_{\alpha}],p_{\alpha}]}{m_{\alpha}}=-
\frac{[p_{\alpha},[y_{\alpha},p_{\alpha}]]}{m_{\alpha}}=-m_{\alpha}y_{\alpha}$$\\
Differentiating above with respect to vector field $x$ we get;
$$[\partial_{x}Iy^{\alpha},p_{\alpha}]+[Iy_{\alpha},x_{\alpha}]=-m_{\alpha}\partial_{x}y^{\alpha}$$
\begin{eqnarray}
&\Rightarrow & \frac{[\partial_{x}Iy^{\alpha},p_{\alpha}]}{m_{\alpha}}+
\frac{[Iy_{\alpha},x_{\alpha}]}{m_{\alpha}}=-\partial_{x}y^{\alpha} \nonumber \\  
&\Rightarrow & I\partial_{x}Iy^{\alpha}=-\partial_{x}y^{\alpha}-
\frac{[Iy_{\alpha},x_{\alpha}]}{m_{\alpha}}. \nonumber
\end{eqnarray}
Using the same argument we have; 
$$I\partial_{y}Ix^{\alpha}=-\partial_{y}x^{\alpha}-
\frac{[Ix_{\alpha},y_{\alpha}]}{m_{\alpha}}.$$ \\
since 
$$[Iy^{\alpha},p_{\alpha}]=- 
\frac{[p_{\alpha},[y_{\alpha},p_{\alpha}]]}{m_{\alpha}}=-m_{\alpha}y_{\alpha}$$\\
now differentiating this with respect to the  vector field $Ix$ we get
$$[\partial_{Ix}Iy^{\alpha},p_{\alpha}]+[Iy_{\alpha},Ix_{\alpha}]=-m_{\alpha}\partial_{Ix}y^{\alpha}$$
\hspace{-2cm}
$$\Rightarrow 
\frac{[\partial_{Ix}Iy^{\alpha},p_{\alpha}]}{m_{\alpha}}+
\frac{[Iy_{\alpha},Ix_{\alpha}]}{m_{\alpha}}=-\partial_{Ix}y^{\alpha}$$\\ 
apply $I$ both sides
\begin{eqnarray}
\hspace{-2cm} &\Rightarrow &
\frac{[[\partial_{Ix}Iy^{\alpha},p_{\alpha}],p_{\alpha}]}{m_{\alpha}^{2}}=-I\partial_{Ix}y^{\alpha} \nonumber \\
\hspace{-2cm}
&\Rightarrow & 
\frac{[p_{\alpha},[\partial_{Ix}Iy^{\alpha},p_{\alpha}]]}{m_{\alpha}^{2}}=I\partial_{Ix}y^{\alpha} \nonumber \\
\hspace{-2cm}
&\Rightarrow & \partial_{Ix}Iy^{\alpha}-
\frac{(p_{\alpha},\partial_{Ix}Iy^{\alpha})p_{\alpha}}{m_{\alpha}^{2}}=I\partial_{Ix}y^{\alpha}. \nonumber
\end{eqnarray}
Since $(p_{\alpha},Iy^{\alpha})=0$ we obtain $(p_{\alpha},\partial_{Ix}Iy^{\alpha})=-(Ix_{\alpha},Iy_{\alpha}).$ Hence 
$$\partial_{Ix}Iy^{\alpha}+
\frac{(Ix_{\alpha},Iy_{\alpha})p_{\alpha}}{m_{\alpha}^{2}}=I\partial_{Ix}y^{\alpha}$$\\
Similarly,
$$\partial_{Iy}Ix^{\alpha}+
\frac{(Iy_{\alpha},Ix_{\alpha})p_{\alpha}}{m_{\alpha}^{2}}=I\partial_{Iy}x^{\alpha}$$
Substituting corresponding expressions in $N_{I}(x,y)$ we get;
$$N_{I}(x,y)=I\partial_{Ix}y^{\alpha}-
\frac{(Ix_{\alpha},Iy_{\alpha})p_{\alpha}}{m_{\alpha}^{2}}-I\partial_{Iy}x^{\alpha}+\frac{(Iy_{\alpha},Ix_{\alpha})p_{\alpha}}{m_{\alpha}^{2}}-\partial_{x}y^{\alpha}+\partial_{y}x^{\alpha}-$$
\hspace{1.5cm}
$$I\partial_{Ix}y^{\alpha}-\partial_{y}x^{\alpha}-
\frac{[Ix_{\alpha},y_{\alpha}]}{m_{\alpha}}+\partial_{x}y^{\alpha}+
\frac{[Iy_{\alpha},x_{\alpha}]}{m_{\alpha}}+I\partial_{Iy}x^{\alpha}=0 \qquad \Box$$
\begin{df} A K\"ahler manifold is a symplectic manifold with an integrable almost complex structure.
\end{df}
\begin{cor}$\mathcal{M}$ is a K\"ahler manifold with K\"ahler structure 
$$\Omega(u,v)=g(u,v)+i\omega (u,v).\quad \Box$$
\end{cor}

\textit{Anknowledgemts} I want to thank A. Klyachko for his patient and valuable ideas and I want to thank Pomona College for their hospitality.

Vehbi E. Paksoy
CMC Mathematics Dept.
Claremont, CA 91711, USA
emrah.paksoy@cmc.edu

\end{document}